\documentclass[11pt]{article}
\usepackage{epsf}
\usepackage{color}
\usepackage{epsfig}
\usepackage{times}

\usepackage{graphicx}
\usepackage{amsmath}
\usepackage{latexsym}
\usepackage{amssymb}
\usepackage{amsfonts}
\usepackage{url}

\usepackage{a4}

\newcommand{\psis}{{\mathcal S}}

\newcommand{\e}{\ensuremath{\mathrm{e}}}

\newcommand{\psib}{\chi}

\begin{document}

\title{Splitting and composition methods with embedded error estimators} 

\author{Sergio Blanes\footnote{Universitat Polit\`ecnica de Val\`encia, Instituto de Matem\'atica Multidisciplinar, 46022 Valencia, Spain. Email: serblaza@imm.upv.es}\,, 
Fernando Casas\footnote{Universitat Jaume I, IMAC, Departament de Matem\`atiques, 12071~Castell\'on, Spain. Email: fernando.casas@uji.es}\,, 
Mechthild Thalhammer\footnote{Leopold--Franzens Universit\"at Innsbruck, Institut f\"ur Mathematik, 6020 Innsbruck, Austria. Email: mechthild.thalhammer@uibk.ac.at}}

\date{}
\maketitle

\begin{abstract}
  We propose new local error estimators for splitting and composition methods. They are based on the construction of lower order schemes
  obtained at each step as a linear combination of the intermediate stages of the integrator, so that the additional computational cost required
  for their evaluation is almost insignificant. These estimators can be subsequently used to adapt the step size along the integration.
  Numerical examples show the efficiency of the procedure. 
  
\end{abstract}

\section{Introduction}

Splitting and composition methods are of particular interest in the numerical integration of differential equations when the vector field is separable into solvable parts or when a low order basic method is known, and the goal is to construct higher order schemes by composing the basic method with fractional time steps \cite{mclachlan02sm,mclachlan06gif}.

Although integrators of this class have a long history in numerical mathematics and have been applied, sometimes with different names, in many different contexts (partial differential equations \cite{strang68otc}, quantum statistical mechanics \cite{takahashi84mco}, chemical physics \cite{feit82sot,gray96sit}, molecular dynamics \cite{verlet67ceo}, celestial mechanics \cite{candy91asi,laskar01hos}, etc.), it has been with the advent of the so-called \emph{Geometric Numerical Integration} that the interest in splitting 
and composition has revived and new and very efficient schemes have been designed in the simulation of physical systems.
The goal in Geometric Numerical Integration is to construct schemes in such a way that the numerical approximation
shares with the exact solution many of its relevant qualitative (very often, geometrical) properties, such as symplecticity, unitarity, orthogonality, etc. \cite{blanes16aci,hairer06gni}. If the basic method possesses (some of) these
geometric properties, so do the schemes obtained by composing them. In addition, when they are used with a constant time step, they show a more favorable error growth
behavior than standard integrators, especially in long term integrations. Symplectic integration schemes for Hamiltonian dynamical systems constitute a classical example
of geometric numerical integrators \cite{sanz-serna94nhp}.

Even in problems where no qualitative properties have to be preserved and/or only short time integrations are required, splitting and composition methods have shown to be 
an excellent option (see e.g. \cite{glowinski16smi} and references therein), even when compared with other standard integrators. 

As is well known, some of the most  popular and efficient standard schemes are embedded methods: the numerical procedure contains, besides the numerical 
approximation $x_n$, a second approximation $\widetilde{x}_n$ (usually of a lower order) obtained from intermediate outputs, so that the difference is used as an estimate of the local error for the less precise result and can subsequently be used for step size control \cite{hairer93sod}. Well known
examples in this area are the class of high order embedded Runge--Kutta methods constructed by Verner \cite{verner78erk} 
(implemented as the DVERK code) and Prince \& Dormand \cite{prince81hoe}, giving rise to the code
DOP853 \cite{hairer93sod}.

Since splitting and composition methods also provide intermediate outputs when computing the numerical approximation at every step, it seems then natural to 
analyze whether these intermediate outputs can also be used along the same lines as standard embedded methods to endow the schemes with a step size control.
We will see that this is indeed the case as long as the splitting scheme involves a sufficiently large number of stages and, furthermore, we will show how to construct explicitly
the lower order approximation $\widetilde{x}_n$ from these intermediate outputs at virtually cost free. 

It is important to remark that, whereas splitting and composition methods implemented with a constant step size are specially well suited in geometric numerical integration
for long time integrations, this is not the case of the variable step size
schemes constructed by applying the strategy proposed here \cite{calvo93tdo}. In any case, the second approximation $\widetilde{x}_n$ is only used to estimate the local error and this
is not propagated along the integration interval.

Of course, the idea of endowing splitting methods with a local error estimator is not new. We can mention in particular references \cite{descombes11ats,descombes16osm},
where a embedded splitting method is constructed for the second-order Strang splitting for stiff evolutionary partial differential equations, and 
\cite{auzinger16asm,auzinger17aho,thalhammer12ans}, where
a controller splitting method of order $r + 1$ is selected and then an integrator of order $r$ is constructed for which a maximal number of compositions coincide with those of the controller. The methods thus built are then applied for the numerical solution of nonlinear parabolic problems with periodic boundary conditions.

By contrast, the approach we follow here allows one, given a splitting or composition method of order $r$, to construct a second,  lower order approximation as a 
linear combination of the outputs generated at the intermediate stages. This is essentially similar to the procedure presented in \cite{blanes06cmf} for computing
cheap approximations to the optimal postprocessor in composition methods with processing, and can be done virtually cost-free. The lower order methods
thus designed can be used to endow some of the most popular splitting and composition schemes with a reliable and easy-to-evaluate error
estimator \cite{blanes02psp,chin97sif,omelyan03sai}

The plan of the paper is the following. In section 2 we briefly summarize the mathematical formalism to be used in the subsequent analysis. 
Then, in section 3 we proceed to obtain estimators for symmetric compositions of second order basic schemes and of a first order method with 
its adjoint, whereas an analogous treatment is discussed in section 4. The relationship between composition and splitting methods, together with
their respective estimators, is treated in section 5. The new estimators are illustrated in section 6 in comparison with other well established techniques.
Finally, section 7 contains some concluding remarks.

\section{Flows and Lie derivatives}

The analysis of splitting and composition methods can be conveniently carried out with the formalism of Lie derivatives. In that case both the exact flow
and the numerical flow corresponding to an integrator, as well as compositions of this integrator, can be associated to the exponential (or
products of exponentials) of operators, just as in the linear case, so that the order conditions can be obtained by applying the familiar 
Baker--Campbell--Hausdorff formula. 

To be more specific, given the initial value problem
\begin{equation}   \label{eq.1.1}
   \dot{x} = f(x), \qquad x_0 = x(0) \in \mathbb{R}^D
\end{equation}
with $f: \mathbb{R}^D \longrightarrow  \mathbb{R}^D$ and flow $\varphi_t$, we can associate with $f$ the
first order differential operator (the Lie derivative)
$L_f$, whose action on differentiable functions $G: \mathbb{R}^d \longrightarrow \mathbb{R}$ is (see \cite[Chap. 8]{arnold89mmo})
\[
   L_f G(x) = \sum_{i=1}^d f_i(x) \frac{\partial G}{\partial x_i},
\]
so that formally 
\begin{equation} \label{nlp4}
  L_{f} = \sum_{i=1}^d f_i \frac{\partial}{\partial x_i}.
\end{equation}
Moreover, one can also introduce an operator $\Phi^t$ acting on  functions $G$ as  \cite{olver93aol}
\begin{equation}\label{nlp2}
 \Phi^t [G](x) = (G \circ \varphi_t)(x). 
\end{equation}
Then, the Taylor series of $G(\varphi_t(x_0))$ at
$t=0$ is given by \cite{hairer06gni,blanes16aci}
\begin{equation}  \label{nlp2a}
  G(\varphi_t(x_0)) = \sum_{k \ge 0} \frac{t^k}{k!} (L_{f}^k G) (x_0) \equiv \exp (t L_{f})[G](x_0),
\end{equation}
and so 
\begin{equation}   \label{exact.flow}
  \Phi^t[G](x) = \exp(t L_{f})[G](x) \equiv \exp(t F)[G](x),  
\end{equation}
where, for the sake of simplicity in the notation, we write $F \equiv L_f$.
If we replace $G$ in (\ref{exact.flow}) by the identity map $\mathrm{Id}(x) = x$, we get for the exact solution of (\ref{eq.1.1})
\begin{equation}  \label{exact.sol}
  \varphi_t(x_0) = \exp 
   (t F)[\mathrm{Id}](x_0).
\end{equation}
In the same way as for the exact flow $\varphi_t$, we can associate to each numerical integrator for a time step $h$, $\chi_h: \mathbb{R}^d \longrightarrow \mathbb{R}^d$,
the operator 
\begin{equation}
\label{eq:X(h)}
X(h) = I + \sum_{n\geq 1} h^n X_n,
\end{equation}
where $I$ denotes the identity operator and each  $X_n$ acts on smooth functions $G$ as
\begin{equation}
\label{eq:X_n}
X_n[G](x) =  \left. \frac{1}{n!}\frac{d^n}{dh^n}\right|_{h=0}  G(\chi_h(x)),
\end{equation}
so that $X(h)[G](x) = (G \circ \chi_h)(x)$. It is then possible to write $X(h)$ formally as the exponential of another operator $Y(h)$, 
\begin{equation}  \label{fint.1}
  (G \circ \chi_h)(x) = X(h)[G](x) = \exp(Y(h))[G](x),
\end{equation}
where
\begin{equation}   \label{fint.0}
Y(h) = \sum_{n\geq 1} h^n Y_n = \log(X(h)). 
\end{equation}
Clearly, the integrator $\chi_h$ is of order $r$ if $\exp(Y(h)) = \exp(h F) $ up to terms $h^r$, or equivalently, if 
\begin{equation*}
 Y_1 = F , \quad \mbox{ and } \quad Y_n = 0 \quad  \mbox{for} \quad 2 \leq n \leq r.
\end{equation*}
Thus, in particular, if $r=1$, then 
\[
  \exp(Y(h)) = \exp \big(h F + h^2 Y_2 + h^3 Y_3 + \mathcal{O}(h^4) \big),
\]
whereas for its \emph{adjoint method}   $\chi^*_h \equiv \chi_{-h}^{-1}$, one has analogously
\[
  (G \circ \chi^*_h)(x) = \exp(-Y(-h))[G](x)
\]  
with
\[
  \exp(-Y(-h)) = \exp \big(h F - h^2 Y_2 + h^3 Y_3 + \mathcal{O}(h^4) \big).
\]
A second-order method $\psis^{[2]}_{h}$ is (time-)symmetric if and only if $(\psis^{[2]}_{h})^* = \psis^{[2]}_{h}$, or equivalently, if 
its corresponding operator has the form
 $Y(h) = h F + h^3 Y_3 +  h^5 Y_5 + \mathcal{O}(h^7)$.

\section{Estimators for composition methods}

\subsection{Composition of symmetric second order methods}

Suppose now that, starting with a basic symmetric second order integrator $\psis^{[2]}_{h}$, we form the composition
\begin{equation} \label{eq:compSS}
  \psi_h = \psis^{[2]}_{h \alpha_s} \circ \cdots \circ \psis^{[2]}_{h \alpha_2} \circ \psis^{[2]}_{h \alpha_1}.  
\end{equation}
If the coefficients $\alpha_1, \ldots, \alpha_s$ satisfy some requirements (the \emph{order conditions}), then $\psi_h$ provides an approximation of order $r$ to the exact solution.
The number of order conditions is considerably reduced for symmetric compositions, i.e.,
\begin{equation}
  \label{eq:sym3}
\alpha_{j}=\alpha_{s-j+1}, \quad \mbox{ for all } \, j
\end{equation}
in (\ref{eq:compSS}). In that case its associated series of differential operators reads
\[
  \Psi(h) = \exp( Y(h \alpha_1)) \,  \exp(Y(h \alpha_2)) \, \cdots \, \exp(Y(h \alpha_2))   \exp(Y(h \alpha_1)), 
\]
where
\begin{equation}  \label{yexp}
    \exp( Y(h \alpha_k)) = \exp \big(h \alpha_k F + h^3 \alpha_k^3 Y_3 + h^5 \alpha_k^5 Y_5 + \mathcal{O}(h^7) \big)
\end{equation}
is the operator associated with $\psis^{[2]}_{h}$. By requiring that
\[
   \Psi(h) = \exp(h F + \mathcal{O}(h^{r+1})),
\]
one gets the order conditions to be satisfied by the coefficients $\alpha_1, \ldots, \alpha_s$ in the composition (\ref{eq:compSS}).  Up to order $r=6$ these conditions read
explicitly 
\begin{eqnarray}  \label{orcon.6}
  &  &  \sum_{j=1}^s \alpha_j = 1, \qquad\quad    \sum_{j=1}^s \alpha_j^3 = 0   \\
  &  &   \sum_{j=1}^s \alpha_j^5 = 0,  \qquad\quad
  \sum_{j=1}^s \alpha_j^3  \left( \left( \sum_{\ell=1}^{j-1} \alpha_{\ell} \right)^2 + \alpha_j  \,  \sum_{\ell=1}^{j-1} \alpha_{\ell}  \right) = 0 . \nonumber
\end{eqnarray}
Notice that, when computing the numerical approximation $x_{n+1} = \psi_h(x_n) \approx x(t_{n+1}) = x(t_n + h)$ with (\ref{eq:compSS}), the procedure also provides
$s-1$ intermediate outputs in addition to $x_n$, i.e., 
\[
  x_{n,k}= \psis^{[2]}_{h \alpha_k} \circ \cdots \circ \psis^{[2]}_{h \alpha_1}x_n, \qquad k=1,\ldots,s-1,
\]
and the question we pose is whether one can obtain another approximation $\widetilde{x}_{n+1}$ of $x(t_{n+1})$ by a linear combination
\begin{equation}   \label{lic1}
   \widetilde{x}_{n+1} = \sum_{k=0}^{s-1} w_k \, x_{n,k} 
\end{equation}
of these intermediate values $x_{n,k}$, with $x_{n,0} = x_n$. It turns out that this is indeed possible, but the highest order of approximation that can be achieved in this way
depends on the number of intermediate stages $s$. The procedure is similar to the technique used in \cite{blanes04otn,blanes06cmf} to construct cheap postprocessors 
for composition methods with processing.

One should note that $w_s$ is not included in the linear combination (\ref{lic1}). Otherwise, only the trivial solution
\[
 w_s=1, \qquad w_k=0, \quad k=0,1,\ldots,s-1
\]
is obtained. 

Our goal is then to find coefficients $w_k$ so that, given a number of stages $s$, the linear combination (\ref{lic1}) is an approximation to $x(t_{n+1})$ of order
$\ell$, or equivalently,
\begin{equation}  \label{lic2}
   w_0 \, I + \sum_{k=1}^{s-1} w_k  \, \prod_{i=1}^k \exp(Y(h \alpha_i)) = \exp(h F) + \mathcal{O}(h^{\ell +1}),
\end{equation}
where $\exp(Y(h \alpha_k))$ is given by (\ref{yexp}) and $\ell$ is as large as possible. 
Since a linear combination of exponential operator is not, in general, a exponential operator, the conditions to be satisfied by $w_k$
can be derived by expanding both terms in (\ref{lic2}) in powers of $h$ and equating their respective coefficients.  Thus, in particular, up to order $\ell = 4$, one has explicitly
\[
  \exp(h F) = I + h F + \frac{h^2}{2} F^2 + \frac{h^3}{3!} F^3  +  \frac{h^4}{4!} F^4 + \mathcal{O}(h^5)
\]
and
\[
\aligned
  & w_0 \, I + \sum_{k=1}^{s-1} w_k  \, \prod_{i=1}^k \exp(Y(h \alpha_i))   =    f_0 I + h f_1 F + \frac{h^2}{2} f_2 F^2  \\
  & \qquad  +h^3 \left( \frac{1}{3!} f_{3,1} F^3 + f_{3,2} Y_3 \right)   +  h^4  \left( \frac{1}{4!} f_{4,1} F^4 + \frac{1}{2} f_{4,2} F \, Y_3 + \frac{1}{2} f_{4,3} Y_3 \, F \right) + \mathcal{O}(h^5),
\endaligned
\]
whence the following system of linear equations results:
\begin{equation}  \label{orcon1}
\aligned
 & f_0 \equiv w_0 + \sum_{k=1}^{s-1} w_k = 1 \\
 & f_1 \equiv \sum_{k=1}^{s-1} w_k  \sum_{j=1}^k \alpha_j = 1 \\
 & f_2 \equiv \sum_{k=1}^{s-1} w_k  \, \Big(\sum_{j=1}^k \alpha_j \Big)^2 = 1 \\
 & f_{3,1} \equiv \sum_{k=1}^{s-1} w_k  \, \Big(\sum_{j=1}^k \alpha_j \Big)^3= 1 \\
 & f_{3,2} \equiv \sum_{k=1}^{s-1} w_k  \, \Big(\sum_{j=1}^k \alpha_j^3 \Big)= 0 \\
 & f_{4,1} \equiv \sum_{k=1}^{s-1} w_k  \, \Big(\sum_{j=1}^k \alpha_j \Big)^4= 1 \\
 & f_{4,2} \equiv \sum_{k=1}^{s-1} w_k  \, \Big( \sum_{j=1}^k \alpha_j^4 + 2 \sum_{j=1}^{k-1} \alpha_j^3 \sum_{\ell = j+1}^k \alpha_{\ell} \Big) = 0 \\
 & f_{4,3} \equiv \sum_{k=1}^{s-1} w_k  \, \Big( \sum_{j=1}^k \alpha_j^4 + 2 \sum_{j=2}^{k} \alpha_j^3 \sum_{\ell = 1}^{j-1} \alpha_{\ell} \Big) = 0 .
\endaligned
\end{equation}
Notice that the first equation is trivially solved in $w_0$, so to achieve an approximation $\widetilde{x}_{n+1}$ of order 4, we have to verify 7 linear equations. More generally,
the total number of equations (in addition to the trivial one) required to achieve a given order $\ell$ is collected in Table \ref{tab.1} for orders $\ell = 1, \ldots 6$.
Strictly speaking, this number is the sum of the dimensions $m_k$, $k \ge 1$, of the subspaces $\mathcal{A}_k$ of the universal enveloping algebra 
$\mathcal{A} = \bigoplus_{k \ge 0} \mathcal{A}_k$ associated to the graded Lie algebra of operators corresponding to the composition method, with $\mathcal{A}_0 = \mbox{span}(I)$
\cite{blanes06cmf}.

\begin{table}  
\begin{center}
  \begin{tabular}{|c||cccccc|} \hline
 Order $\ell$ & 1 & 2 & 3 & 4 & 5 & 6  \\ \hline\hline 
$\mathcal{SS}$ & 1 & 2 & 4 & 7 & 12 &  20   \\ \hline
method-adjoint &  1  &  3  &  7  &  15  &  31  &  63 \\ \hline
splitting  &  2  &  6  &  14  &  30  &  62  &  126 \\ \hline
\end{tabular}
\end{center}
\caption{Number of order conditions, in additional to the trivial one for $w_0$, required by a linear combination of intermediate outputs to achieve order $\ell$ for symmetric
compositions of 2nd-order symmetric schemes ($\mathcal{SS}$), compositions of a first order method with its adjoint (\ref{eq:compint}) (method-adjoint) and a splitting method
(\ref{split1}) (splitting).}
\label{tab.1}
\end{table}

Next we analyze in detail the construction of numerical schemes of orders 3, 4 and 5 within this approach to be used as error estimators for symmetric
compositions of the form (\ref{eq:compSS}). 

\paragraph{Third-order estimators.} 
Only the first five equations in (\ref{orcon1}) have to be satisfied to get order three. This can be achieved if the composition (\ref{eq:compSS}) has at least $s=5$. For  $s=5$, when the
symmetry of the coefficients (\ref{eq:sym3}) (i.e., $\alpha_5 = \alpha_1$, $\alpha_4 = \alpha_2$) and the order conditions of a 4th-order composition 
(i.e., equations in the first line of (\ref{orcon.6})) are taken into account, then the unique solution of the system is given by
\begin{equation}  \label{s5}
\aligned
&    w_1 = w_4 = \frac{g_2 (1-g_2)}{g_1(g_1 - 1)- g_2(g_2 - 1)}, \qquad
     w_2 = w_3 = 1-w_1
			\\
     & g_1 = \alpha_1, \qquad g_2 = \alpha_1 + \alpha_2
\endaligned     
\end{equation}
so that $w_0 = -1$. A popular (and efficient) 4th-order composition method within this class is the one devised by Suzuki
 \cite{suzuki91gto}, with coefficients
\begin{equation}   \label{suz1}
  \alpha_1 = \alpha_2  =   \frac{1}{4-4^{1/3}}, \qquad \alpha_3 = \frac{1}{1-4^{2/3}}, 
\end{equation}
so that its third-order estimator reads
\begin{equation}   \label{suz2}
  \widetilde{x}_{n+1} = - x_n + w_1 (x_{n,1} + x_{n,4}) + w_2 (x_{n,2} + x_{n,3}).
\end{equation}  
Another widely used 4th-order method involving $s=7$ stages is due to McLachlan \cite{mclachlan02foh}, with coefficients
\[
  \alpha_1 = \alpha_2 = \alpha_3 = \frac{1}{6-6^{1/3}}, \qquad \alpha_4 = \frac{1}{1-6^{2/3}}.
\]
Its corresponding estimator now involves a free parameter, which can be taken to be $w_3$, and reads
\[
   \widetilde{x}_{n+1} = - x_n + w_1 (x_{n,1} + x_{n,6}) + w_2 (x_{n,2} + x_{n,5}) + w_3 (x_{n,3} + x_{n,4}).
\]  
Here
\[
\aligned
  & w_1 =  \frac{g_2 (1-g_2) + w_3 (g_2(g_2 -1) - g_3(g_3-1))}{g_1(g_1 - 1)- g_2(g_2 - 1)} \\
   &   w_2 = \frac{g_1 (g_1 - 1) - w_3 (g_1(g_1 -1) + g_3(g_3-1))}{g_1(g_1 - 1)- g_2(g_2 - 1)}, 
\endaligned
\]
with $g_1 = \alpha_1$, $g_i = g_{i-1} + \alpha_i$, $i=2,3$.   

The same strategy can also be applied to the  popular 4th-order 3-stage Yoshida's method \cite{yoshida90coh}
\begin{equation}   \label{eq:Yoshida3}
  \phi^{[4]}_h = 
	\psis^{[2]}_{h \alpha_1} \circ \psis^{[2]}_{h \alpha_2} \circ \psis^{[2]}_{h \alpha_1}, 
\end{equation}
with
\[
\displaystyle
  \alpha_1= \frac{1+2^{-5/3}}{2+2^{1/3}+2^{-1/3}} \pm
	\frac{i}4 \, \frac{\sqrt{3}}{1+2^{2/3}+2^{-2/3}}, \qquad
  \alpha_2=1-2\alpha_1
\]
which is known to lead to small errors when complex coefficients are taken \cite{blanes13oho}. 
Since  only three intermediate outputs per step are available, one needs at least two steps of it as if it were one single method, i.e., one can take as 
integrator the composition
\begin{equation}\label{eq:Suzuki4}
  \phi^{[4]}_h = 
	\psis^{[2]}_{h \alpha_1/2} \circ \psis^{[2]}_{h \alpha_2/2} \circ 
	\psis^{[2]}_{h \alpha_1/2}  \circ 
	\psis^{[2]}_{h \alpha_1/2} \circ \psis^{[2]}_{h \alpha_2/2} \circ 
	\psis^{[2]}_{h \alpha_1/2}. 
\end{equation}
In this case the corresponding estimator reads
\[
   \widetilde{x}_{n+1} = - x_n + w_1 (x_{n,1} + x_{n,5}) + w_2 (x_{n,2} + x_{n,4}), 
\]  
with
\[ 
   w_1 = \frac{1- \alpha_1^2}{\alpha_2}, \qquad w_2 = 1-w_1.
\]   
We can adopt the terminology of embedded Runge--Kutta methods \cite{hairer93sod} and denote the previous compositions with their 
respective estimators as methods of order 4(3).

\paragraph{Compositions of order 6(4).}
 
To get linear combinations (\ref{lic1}) of order four one has to solve the whole set of equations (\ref{orcon1}). Although in principle this would require $s=8$ stages,
it turns out that if the underlying
time-symmetric composition (\ref{eq:compSS}) satisfies the order conditions up to order 6 given by (\ref{orcon.6}) 
with the minimum number of stages ($s=7$), one gets a unique
solution of the form
\[
   \widetilde{x}_{n+1} =  x_n + w_1 (x_{n,1} - x_{n,6}) + w_2 (x_{n,2} - x_{n,5}) + w_3 (x_{n,3} - x_{n,4}), 
\]  
where $w_i$ can be expressed analytically in terms of the $\alpha_i$ coefficients of the composition. For the particular method found by Yoshida
\cite{yoshida90coh}, with coefficients
\[
\aligned
 & \alpha_1 = 0.78451361047755726382, \quad \;\;\;
   \alpha_2 = 0.23557321335935813369 \\
 & \alpha_3 =-1.17767998417887100695, \quad
   \alpha_4 = 1-2(\alpha_1+\alpha_2+\alpha_3) 
\endaligned
\]
one has
\[
\aligned
 & w_1 =-0.90983233007647709242, \\
 & w_2 = 2.16331188722978237305, \\
  & w_3 = 0.55695580387159066608.
\endaligned  
\]
The same strategy can be applied of course if 6th-order compositions with more stages are considered. For instance, we have found an estimator within this class
for the symmetric method proposed by Kahan \& Li \cite{kahan97ccf}, with $s=9$ stages.  
  
\paragraph{Compositions of order 6(5) and 8(5).}

A system of 13 linear equations has to be solved for getting an estimator of order five. Although not all of them are independent when the time-symmetry
and the order conditions for the underlying composition are introduced, at least $s=11$ stages are necessary. Starting from the 6th-order symmetric composition 
obtained by Sofroniou \& Spaletta \cite{sofroniou05dos} with coefficients

\begin{equation}   \label{m65al}
\aligned
 & \alpha_1 = 0.21375583945878254555, \quad
   \alpha_2 = 0.18329381407425713911 \\
 & \alpha_3 = 0.17692819473098943795, \quad
   \alpha_4 =-0.44329082681170215849  \\
 & \alpha_5 = 0.11728560432865935385, \quad
   \alpha_6 = 1-2(\alpha_1+\alpha_2+\alpha_3+\alpha_4+\alpha_5), 
\endaligned
\end{equation}  
there is just one set of coefficients satisfying all the order conditions. The resulting
method of order 6(5) is of the form
\begin{equation}   \label{m65}
    \widetilde{x}_{n+1} =  -x_n + \sum_{i=1}^5 w_i  \, (x_{n,i} + x_{n,11-i}),  
\end{equation}  
with
\[
\begin{array}{lcl}
  w_1 = -4.70925883588386976399  &  & 
   w_2 = 24.61043285614692442695 \\
   w_3 =-19.39218824966918044634 &  & 
 w_4 =  6.17441462307605721006 \\
   w_5 = -5.68340039366993142668  &  &
\end{array}
\]  
The same strategy can be applied to compositions (\ref{eq:compSS}) of order 8. A well known example within this class is the symmetric method
proposed by Kahan \& Li  \cite{kahan97ccf}
with $s=17$ and coefficients 
\begin{equation}   \label{krcl0}
\begin{array}{lcl}
 \alpha_1 = 0.13020248308889008088, &  & 
   \alpha_2 = 0.56116298177510838456 \\
 \alpha_3 =-0.38947496264484728641,  &  & 
   \alpha_4 = 0.15884190655515560090 \\
  \alpha_5 =-0.39590389413323757734, &  & 
   \alpha_6 = 0.18453964097831570709 \\
  \alpha_7 = 0.25837438768632204729,  &  & 
   \alpha_8 = 0.29501172360931029887 \\
  \alpha_9 = 1-2(\alpha_1+\cdots+\alpha_8), &  &  
\end{array}
\end{equation}  
the estimator reads
\begin{equation}   \label{krcl1}
   \widetilde{x}_{n+1} =  -x_n + \sum_{i=1}^8 w_i  \, (x_{n,i} + x_{n,17-i}),  
\end{equation}  
with
\[
\begin{array}{lcl}
  w_1 =-2.77811433347582461058,  &  & 
  w_2 = 1.43336350604816157334\\
    w_3 =-2.35490307436226712937,  & & 
   w_4 = 0.27249477875971647996 \\
   w_5 = 3.09204406313073660493,  & & 
  w_6 = 1.33511505989947708172 \\
   w_7 = 0, & & 
   w_8 = 0.
\end{array}
\]  
The  DOP853 algorithm based on a 12-stage RK8(6) method by Dormand \& Prince (announced but not published in \cite{dormand89prk}),
where the embedded 6th-order method is replaced by a pair of embedded methods of order five and three by Hairer \& Wanner \cite{hairer93sod}), is one of  the most efficient schemes within this framework. In comparison, the previous composition method involves more stages, but on the other hand does
 not require to keep up to 12 vectors in memory.

As a matter of fact, we can apply the same strategy to the 8th-order composition method considered here and construct a second estimator of
order 3 to avoid any possible over-estimation of the error. One possible 3th-order estimator is given by 
\begin{equation}   \label{krcl2}
  \widetilde{x}_{n+1}^{[3]}= -x_n + w_1 (x_{n,1} + x_{n,16}) + w_7 (x_{n,7} + x_{n,10}),
\end{equation}  
with $w_1$, $w_7$ verifying
\begin{eqnarray*}
   &  &  w_1 + w_7 = 1 \\  
   &  & g_1(g_1 -1) w_1 + g_7 (g_7 -1) w_7 = 0
\end{eqnarray*}
where $g_1 = \alpha_1,  \  g_7 = \alpha_1 + \cdots + \alpha_7 $,
i.e. 
\[
w_1=1.828514038642564624 , \qquad  w_7=-0.828514038642564624.
\]
We then have two error estimators for the scheme (\ref{eq:compSS}) with coefficients (\ref{krcl0}),
\[
  \mbox{\textit{err}}_5 = \| \widetilde{x}_{n}^{[5]}-x_n\|={\cal O}(h^6), \qquad
  \mbox{\textit{err}}_3 = \| \widetilde{x}_{n}^{[3]}-x_n\|={\cal O}(h^4).
\]
Applying now the same strategy as in  \cite{hairer93sod}, we consider
\[
  \mbox{\textit{err}} = \mbox{\textit{err}}_5 \cdot \frac{\mbox{\textit{err}}_5}{\sqrt{\mbox{\textit{err}}_5^2+0.01\cdot \mbox{\textit{err}}_3^2}}  ={\cal O}(h^8)
\]
as an error estimator that behaves asymptotically like the global error of the method. 

Notice that  we can obtain error estimators for other composition schemes in a similar way. For example, at order eight one can find in the literature methods with up to 21 stages \cite{hairer06gni,kahan97ccf,sofroniou05dos}, and their relative performance depend on the particular problem to solve as well as on the symmetric second order scheme used as the basic scheme for the composition.

\subsection{Composition of a first order method with its adjoint}

Higher order methods can also be obtained by composing a first order basic method $\chi_h$ and its adjoint $\chi^*_h=\chi_{-h}^{-1}$,
\begin{equation} \label{eq:compint}
\psi_h =  \chi_{\alpha_{2s} h} \circ \chi^*_{\alpha_{2s-1}h}\circ  \cdots  \circ
\chi_{\alpha_{2}h}  \circ  \chi^*_{\alpha_{1}h},
\end{equation}
with appropriately chosen real coefficients $(\alpha_1,\ldots,\alpha_{2s})$. The associated series of differential
operators is of the form
\begin{equation}  \label{ser.comp1}
  \Psi(h)  =   \e^{-Y(-h\alpha_1)} \, \e^{Y(h \alpha_2)} \cdots \, \e^{-Y(-h\alpha_{2s-1})} \, \e^{Y(h \alpha_{2s})},
\end{equation}
where $Y(h \alpha_k) = h \alpha_k F + h^2 \alpha_k^2 Y_2 + h^3 \alpha_k^3 Y_3 + \mathcal{O}(h^4)$.
Again, by requiring that $\Psi(h) = \exp(h F + \mathcal{O}(h^{r+1}))$, one gets the order conditions to be satisfied by the coefficients to achieve order $r$.
These order conditions are considerably simplified if $\alpha_{2s-j+1} = \alpha_j$ for all $j$. In that case the composition (\ref{eq:compint}) is time-symmetric. 

As with symmetric compositions of symmetric second order schemes, here we can also take a linear combination 
\begin{equation}  \label{lcmad}
  \widetilde{x}_{n+1} = w_0 \, x_n + \sum_{k=1}^{2s-1} w_k \, x_{n,k}
\end{equation} 
of intermediate outputs
\[
   x_{n,2i-1} = \chi_{\alpha_{2i-1}}^*(x_{n,2i-2}), \qquad  x_{n,2i} = \chi_{\alpha_{2i}}(x_{n,2i-1}),
\]
to produce an approximation of order $\ell < r$ to be used as an error estimator for the composition (\ref{eq:compint}). The coefficients $w_k$ can
be determined by requiring that
\[
   w_0 I + w_1 \, \e^{-Y(-h \alpha_1)} + w_2 \, \e^{-Y(-h \alpha_1)} \e^{Y(h \alpha_2)} + \cdots = \exp(h F) + \mathcal{O}(h^{\ell + 1}).
\]   
By expanding the product of exponentials we get the number of conditions the $w_k$ have to satisfy at a given order in a similar way as with compositions of 2nd-order
symmetric methods. This number is collected in Table \ref{tab.1}.

In particular, 8 linear equations are required to get a 3rd-order approximation in this way. Since several efficient 4th-order methods of this class  with up to
6 stages (or 12 intermediate outputs) are available in the literature, it is in principle possible to get third order estimators for them (even with free parameters for
optimization). As an illustration, for
the symmetric 4th-order method (\ref{eq:compint}) with $s=6$ and coefficients
\begin{equation}  \label{MetAdj1}
\begin{array}{lcl}
 \alpha_1 =  0.08298440641740484666, & &  \alpha_2 = 0.16231455076686615333 \\
 \alpha_3 =  0.23399525073150184666, &  &  \alpha_4 = 0.37087741497957699562 \\
 \alpha_5 = -0.40993371990192559562, &  &  \alpha_6 = 0.05976209700657575333
\end{array}
\end{equation} 
we propose the linear combination (\ref{lcmad}) with $w_0 = -1$ and
\begin{equation}  \label{MetAdj2}
\begin{array}{lcl}
   w_1 = 1.48889386198802799037, &  &  w_2 = -0.03049911761922725390 \\
   w_3 =-0.32603028933442750875, &  &  w_4 = -0.05468276894167474320 \\
   w_5 =-0.02746220037522580999, &  &  w_6 = -0.10043897143494534902 \\
   w_{12-i} = w_i, \quad i=1,\ldots,5. &  & 
\end{array}
\end{equation}

\section{Estimators for splitting methods}

If $f$ in equation (\ref{eq.1.1}) can be \emph{split} as $f = \sum_{i=1}^m f^{[i]}$ for certain functions
$f^{[i]}: \mathbb{R}^D \longrightarrow  \mathbb{R}^D$, in such a way that the equations
\begin{equation}   \label{eq.1.2}
   \dot{x} = f^{[i]}(x), \qquad x_0 = x(0) \in \mathbb{R}^D, \qquad i=1, \ldots, m
\end{equation}
can be integrated exactly, with solutions $x(h) = \varphi_h^{[i]}(x_0)$ at $t = h$, then the basic first-order
method in the composition (\ref{eq:compint}) can be taken simply as
\begin{equation}   \label{eq.1.3}
   \psib_h = \varphi_h^{[m]} \circ \cdots \circ \varphi_h^{[2]} \circ \varphi_h^{[1]},
\end{equation}
whereas its adjoint is just the reversed composition
\begin{equation}   \label{eq.1.3a}
   \psib_h^* = \varphi_h^{[1]} \circ \varphi_h^{[2]} \circ \cdots \circ \varphi_h^{[m]}.
\end{equation}
For $m=2$, i.e., when $f(x)$ is decomposed in just two pieces,
\[
   f=f^{[1]}+f^{[2]},   
\]
one could also consider a time-symmetric composition 
\begin{equation}   \label{split1}
  \psi_h = \varphi_{b_{s+1} h}^{[2]} \circ \varphi_{a_{s} h}^{[1]} \circ \varphi_{b_{s} h}^{[2]}  \circ \cdots \varphi_{b_{2} h}^{[2]} \circ 
     \varphi_{a_{1} h}^{[1]}  \circ \varphi_{b_{1} h}^{[2]} 
\end{equation}
with appropriately chosen coefficients $a_i$, $b_i$ verifying
\[
   a_{s+1-j} = a_j, \qquad b_{s+2-j} = b_j, \qquad j=1,2,\ldots
\]   
to achieve a prescribed order. Here it is also possible to take advantage of the 
intermediate outputs to construct a lower order approximation which may be used as an error estimator for the integrator (\ref{split1}). In this case it has the form
\begin{equation}  \label{lcsplit}
  \widetilde{x}_{n+1} = w_0 \, x_n + \sum_{k=1}^{2s} w_k \, x_{n,k},
\end{equation} 
with
\[
   x_{n,2i-1} =  \varphi_{b_i h}^{[2]}(x_{n,2i-2}), \qquad  x_{n,2i} =\varphi_{a_i h}^{[1]}(x_{n,2i-1}).
\]
As before,
the analysis can be carried out with the associated series of differential operators, which in this case reads
\[
  \Psi(h) = \exp(b_1 h B) \, \exp(a_1 h A) \cdots \exp(b_s h B) \, \exp(a_s h A) \, \exp(b_{s+1} h B),
\]
where $A$ and $B$ denote the Lie derivatives corresponding to $f^{[1]}$ and $f^{[2]}$, respectively:
\[
   A \equiv \sum_{i=1}^D f_i^{[1]}(x) \frac{\partial}{\partial x_i}, \qquad\quad    B \equiv \sum_{i=1}^D f_i^{[2]}(x) \frac{\partial}{\partial x_i}.
\]
Analogously, the conditions to be satisfied by the $w_i$ are determined by expanding the exponentials in
\[
  w_0 I + w_1 \e^{b_1 h B} + w_2 \e^{b_1 h B} \e^{a_1 h A}  + \cdots = \exp(h F) + \mathcal{O}(h^{\ell +1}).
\]  
The number to achieve a given order is collected in Table \ref{tab.1} (last line).

Now  a system of 15 equations have to be satisfied by the coefficients $w_i$ in the linear combination (\ref{lcsplit}) to achieve order 3. As in the preceding cases,
we can take several efficient splitting methods of the form (\ref{split1}) involving enough intermediate steps and construct estimators for them. In particular, for the 4th-order symmetric splitting scheme designed by Blanes \& Moan \cite{blanes02psp},
with 12 intermediate outputs
\begin{equation}   \label{prk1}
  \psi_h = \varphi_{b_1 h}^{[2]} \circ \varphi_{a_1 h}^{[1]} \circ \cdots \varphi_{a_3 h}^{[1]} \circ \varphi_{b_4 h}^{[2]} \circ \varphi_{a_3 h}^{[1]} \cdots  \varphi_{a_1 h}^{[1]}
     \circ \varphi_{b_1 h}^{[2]}
\end{equation}
and coefficients
\[
\begin{array}{lcl}
    b_1 = 0.07920369643119565, &  &  a_1 = 0.209515106613361 \\
    b_2 = 0.35317290604977372, &  &  a_2 =-0.143851773179818 \\
   b_3 =-0.04206508035771952,  &  &  a_3 = 1/2-(a_1+a_2) \\  
    b_4 = 1-2(b_1+b_2+b_3)   &  & 
\end{array}
\]
we propose the linear combination
\begin{equation}   \label{prk2}
  \widetilde{x}_{n+1,k} = -x_{n,0} + \sum_{i=1}^5 w_i (x_{n,i} + x_{n,13-i})
\end{equation}
solving all order conditions with
\begin{equation}  \label{est.1}
\begin{array}{lclcl}
  w_1 = 1, &  &      w_2 = 0.43458657385433203071,  &   &  \\
  w_3 = - w_2, &  &   w_4 = 0.27273581001405423884, &  &
  w_5 = -w_4.
\end{array}
\end{equation}
Another particularly efficient 4th-order splitting method designed for systems of the form
\begin{equation}   \label{rkn1}
   \ddot{y} = g(y),  \qquad y \in \mathbb{R}^D
\end{equation}
when written as a first order system
\[
\frac{d}{dt} \left(
\begin{array}{c}
 y \\  \dot y
\end{array} \right)=
\underbrace{\left(
\begin{array}{c}
 \dot y \\ 0
\end{array} \right)}_{f^{[1]}} +
\underbrace{\left(
\begin{array}{c}
 0 \\ g(y)
\end{array} \right)}_{f^{[2]}} 
\]
corresponds to the composition (\ref{prk1}) with
\begin{equation}   \label{rkn3}
\begin{array}{lcl}
%
	 b_1 = 0.082984406417404, &  &  a_1 = 0.245298957184271 \\
     b_2 = 0.396309801498368, &  &  a_2 = 0.604872665711078 \\
    b_3 =-0.039056304922348, &  &  a_3 = 1/2-(a_1+a_2) 
	\\
    b_4 =  1-2(b_1+b_2+b_3). &  &  
\end{array}
\end{equation}
In this case the estimator has also the form (\ref{prk2}) with
\begin{equation}   \label{est.2}
\begin{array}{lclcl}
  w_1 = 1, &  &      w_2 = 0.43541552923952936004,  &  &  \\
  w_3 = - w_2, &  &   w_4 =-0.17978889668391821731, &  & 
   w_5 = -w_4.
\end{array}
\end{equation}
This splitting method, as well as the error estimator, can also be used to integrate in time the Schr\"odinger equation
\[
 i \frac{\partial}{\partial t}\psi = \left(-\frac{1}{2m} \Delta + V(x) \right) \psi,
\]
where $m$ is the reduced mass, $\Delta$ is the Laplacian operator
and $V(x)$ is the potential. After spatial discretisation one has to solve a linear system of ODEs
\[
 i \dot u = (A + B)u, \qquad u_0\in \mathbb{C}^D,
\]
where $A$ corresponds to the spatial discretization of the kinetic part and $B$ to the potential part. Here $B$ is a diagonal matrix in the coordinates space, whereas $A$ is diagonal in the momentum space, so fast Fourier transform
(FFT) algorithms ${\cal F}$ can be used to compute the action of a $A$ on a vector, $Au=F^{-1}D_AFu$, with $D_A$ a diagonal matrix.

\section{Connection between splitting and composition}

Splitting and composition methods for system $\dot{x} = f^{[1]}(x) + f^{[2]}(x)$ are closely connected. On the one hand, if $\psis^{[2]}_{h} = \varphi_{h/2}^{[2]} \circ \varphi_h^{[1]}
\circ \varphi_{h/2}^{[2]}$ or $\psis^{[2]}_{h} = \varphi_{h/2}^{[1]} \circ \varphi_h^{[2]} \circ \varphi_{h/2}^{[1]}$, then the composition scheme (\ref{eq:compSS}) can be written as
(\ref{split1}), although the opposite is not true in general. On the other hand, if 
$\chi_h =  \varphi^{[2]}_{h}\circ \varphi^{[1]}_{h}$, then $\chi_h^* = \varphi^{[1]}_{h} \circ \varphi^{[2]}_{h}$ and the composition 
(\ref{eq:compint}) reads
\begin{equation}
  \label{eq:MetAdj}
\psi_h = 
\big(\varphi^{[2]}_{\alpha_{2s} h}\circ 
\varphi^{[1]}_{\alpha_{2s}h}\big) \circ 
\big(\varphi^{[1]}_{\alpha_{2s-1} h}\circ 
\varphi^{[2]}_{\alpha_{2s-1}h}\big) \circ 
 \cdots\circ
\big(\varphi^{[2]}_{\alpha_{2} h}\circ 
\varphi^{[1]}_{\alpha_{2}h}\big) \circ 
\big(\varphi^{[1]}_{\alpha_{1} h}\circ 
\varphi^{[2]}_{\alpha_{1}h}\big) .
\end{equation}
Since $\varphi^{[i]}_{h}$ ($i=1,2$) are exact flows, then they verify\footnote{This property is not satisfied, in general, if the exact flows are replaced by numerical approximations.}
$
 \varphi^{[i]}_{\beta h} \circ\varphi^{[i]}_{\delta h}= 
 \varphi^{[i]}_{(\beta+\delta) h}, 
$
and the method can be rewritten as the splitting scheme 
\begin{equation}
  \label{eq:splitting}
\psi_h = \varphi^{[2]}_{b_{s+1} h}\circ
 \varphi^{[1]}_{a_{s}h}\circ \varphi^{[2]}_{b_{s} h}\circ
 \cdots\circ
 \varphi^{[2]}_{b_{2}h}\circ
 \varphi^{[1]}_{a_{1}h} \circ \varphi^{[2]}_{b_{1}h},
\end{equation}
if $b_{1}=\alpha_{1}$ and 
\begin{equation}
  \label{eq:abalpha}
a_j = \alpha_{2j} + \alpha_{2j-1}, \qquad\quad  
b_{j+1} = \alpha_{2j+1} + \alpha_{2j}, \qquad j=1,\ldots,s
\end{equation}
 (with $\alpha_{2s+1}=0$). Conversely, any integrator of the form (\ref{eq:splitting}) with 
 $\sum_{i=1}^{s} a_i = \sum_{i=1}^{s+1} b_i$ can be expressed in the form (\ref{eq:compint}) 
 with $\chi_h =  \varphi^{[2]}_{h}\circ \varphi^{[1]}_{h}$ and 
\begin{equation}
  \label{eq:alphaab}
\begin{array}{l}
\alpha_{2s}=b_{s+1}, \\
\alpha_{2j-1}=a_{j}-\alpha_{2j}, \qquad
\alpha_{2j-2}=b_{j}-\alpha_{2j-1}, \qquad j=s,s-1,\ldots,1, \nonumber
\end{array}
\end{equation}
with $\alpha_{0}=0$ for consistency.
Nevertheless, the intermediate outputs are different in each implementation as well as the number of order
conditions for the estimators. In general this number grows faster with the order for splitting methods. Moreover,
implementing the splitting scheme $\psi_h$ as a composition method is in general more costly because
explicitly obtaining the intermediate values requires the computation of additional basic flows. In more detail,
suppose we write \eqref{eq:splitting} as a composition: 
\begin{eqnarray*}  
\psi_h &=&  \cdots\circ
 \varphi^{[2]}_{b_{2}h}\circ
 \varphi^{[1]}_{a_{1}h} \circ \varphi^{[2]}_{b_{1}h} \\
 &=&  \cdots\circ
 \varphi^{[2]}_{(b_{2}-(a_{1}-b_1))h}\circ
 \underbrace{\varphi^{[2]}_{(a_{1}-b_1)h}\circ
 \varphi^{[1]}_{(a_{1}-b_1)h}}_{\chi_{(a_{1}-b_1)h}}  \circ
 \underbrace{\varphi^{[1]}_{b_{1}h} \circ \varphi^{[2]}_{b_{1}h}}_{\chi^*_{b_1h}} .
\end{eqnarray*}
Then, for the first intermediate output we have
\[
  x_{n+1,1}= \chi^*_{b_1h}(x_{n,0})=\varphi^{[1]}_{b_{1}h} \circ \varphi^{[2]}_{b_{1}h}(x_{n,0}).
\]
However, whereas obviously $\varphi^{[1]}_{(a_{1}-b_1)h} \circ \varphi^{[1]}_{b_{1}h}=\varphi^{[1]}_{a_{1}h}$, the computational cost of computing 
$z=\varphi^{[1]}_{b_{1}h}(y)$ and then $\varphi^{[1]}_{(a_{1}-b_1)h}(z)$ can be in many cases up to twice more costly than directly evaluating $\varphi^{[1]}_{a_{1}h}(y)$. 

For example, taking this composition for solving the Schr\"odinger equation requires the computation of $s$ additional inverse FFTs with respect to the same scheme written as a splitting method. Similarly, taking a composition with the symmetric second order scheme $\psis^{[2]}_{h} = \varphi_{h/2}^{[2]} \circ \varphi_h^{[1]}\circ \varphi_{h/2}^{[2]}$ requires the same number of FFTs as the corresponding splitting composition, but 
taking instead
$\psis^{[2]}_{h} = \varphi_{h/2}^{[1]} \circ \varphi_h^{[2]} \circ \varphi_{h/2}^{[1]}$ as the basic scheme, requires $s$ additional inverse FFTs for the intermediate outputs because $\varphi_{h}^{[1]}$ carries the costly part of the scheme.

A noteworthy exception is the case in which $f^{[1]}$ and $f^{[2]}$ originate from a partitioned ordinary differential equation of the form
\begin{equation}  \label{eq:GpFq}
  \dot{q}=g(p), \qquad \dot{p}=f(q).
\end{equation}
The system can then be written as 
\[
\frac{d}{dt} \left(
\begin{array}{c}
 q \\p
\end{array} \right)=
\underbrace{\left(
\begin{array}{c}
 g(p) \\ 0
\end{array} \right)}_{f^{[1]}} +
\underbrace{\left(
\begin{array}{c}
 0 \\ f(q)
\end{array} \right)}_{f^{[2]}} 
\]
and
\[
 \varphi^{[1]}_{b_{1}h} \left(
\begin{array}{c}
 q_n \\p_n
\end{array} \right) =\left(
\begin{array}{c}
 q_n + b_{1}h g(p_n)\\p_n
\end{array} \right) ,\
 \varphi^{[1]}_{(a_{1}-b_1)h} \left(
\begin{array}{c}
 q_n \\p_n
\end{array} \right) =\left(
\begin{array}{c}
 q_n + (a_{1}-b_1)h g(p_n)\\p_n
\end{array} \right) ,\qquad
\]
where the same evaluation $g(p_n)$ is used in both cases.

The algorithm corresponding to the splitting method (\ref{eq:splitting}) 
for the step $(q_0,p_0) \mapsto (q_1,p_1)$
reads 
\begin{eqnarray*}
  { Q}_0  &  =  &  q_0, \;\;\;\;\;
  { P}_0 = p_0  \nonumber  \\
     &  \mbox{for}   &  i=1, \dots, s    \nonumber  \\
      &   &  { Q}_{2i-1} = { Q}_{2i-2}    \\
      &   &  { P}_{2i-1} = { P}_{2i-2} + h b_i
            f({ Q}_{2i-1})    \nonumber  \\
      &   &  { Q}_{2i} = { Q}_{2i-1} + h a_i
            g({ P}_{2i-1})    \\
      &   &  { P}_{2i} = { P}_{2i-1}   \nonumber  \\
  q_1  &  =  &  { Q}_{2s}, \;\;\;\;\;
  p_1 = { P}_{2s} + h b_{s+1} f({ Q}_{2s}),  \nonumber
\end{eqnarray*}
so that it can be seen as an explicit partitioned Runge--Kutta method. On the other hand,
the composition \eqref{eq:compint} with \eqref{eq:MetAdj} leads to the algorithm
\begin{eqnarray*}
  { Q}_0  &  =  &  q_0, \;\;\;\;\;
  { P}_0 = p_0  \nonumber  \\
     &  \mbox{for}   &  i=1, \dots, s    \nonumber  \\
      &   &  { P}_{2i-1} = { P}_{2i-2} + h \alpha_{2i-1}
            f({ Q}_{2i-2})    \nonumber  \\
      &   &  { Q}_{2i-1} = { Q}_{2i-2}  + h \alpha_{2i-1}
            g({ P}_{2i-1})     \\
      &   &  { Q}_{2i} = { Q}_{2i-1} + h \alpha_{2i}
            g({ P}_{2i-1})     \\
      &   &  { P}_{2i} = { P}_{2i-1}  + h \alpha_{2i}
            f({ Q}_{2i})   \nonumber  \\
  q_1  &  =  &  { Q}_{2s}, \;\;\;\;\;
  p_1 = { P}_{2s}   \nonumber
\end{eqnarray*}
requiring exactly the same evaluations of $f$ and $g$. If in addition $g(p)=p$ (i.e., if we are solving the second order differential equation $\ddot{q}=f(q)$), then 
the estimator for \eqref{eq:GpFq} takes the form (for appropriate choices of the parameters $w_i$)
\begin{eqnarray*}  
  \left(
\begin{array}{c}
\tilde q_{n+1} \\\tilde p_{n+1}
\end{array} \right) &=&  
  w_0 \left(
\begin{array}{l}
 q_n \\	p_n 
\end{array} \right) + 
\sum_{i=1}^{s}\left(
w_{2i-1}\left(
\begin{array}{c}
 Q_{2i-1} \\P_{2i-1}
\end{array} \right) +
w_{2i}
  \left(
\begin{array}{c}
 Q_{2i} \\P_{2i}
\end{array} \right) \right)  \\
&  & +w_{2s-1}\left(
\begin{array}{c}
 Q_{2s-1} \\P_{2s-1}
\end{array} \right) 
 =
  \left(
\begin{array}{l}
\displaystyle  q_n + h p_n + h^2 \sum_{i=1}^{s}\delta_{i}f(Q_{2i-2})\\
\displaystyle 	p_n + h \sum_{i=1}^{s}\gamma_{i}f(Q_{2i-2})
\end{array} \right)
\end{eqnarray*}
in a similar way as for embedded Runge--Kutta--Nystr\"om methods.
In any case, other choices of $\delta_{i},\gamma_{i}$ can also lead to estimators associated to a given $s$-stage composition scheme \cite{calvo93tdo}, and that can not be obtained by taking intermediate outputs.

\section{Numerical examples}

In this section we analyze the accuracy and reliability of the estimators presented in this work in comparison with other well established schemes
for a simple example. Specifically, the methods (and notation) we consider are the following:
\begin{itemize}
	\item $\mathcal{RKN}_{6}43$: The 6-stage 4th-order splitting method (\ref{prk1}) for systems of the form (\ref{rkn1}) with the 
           3rd-order estimator (\ref{est.2}).
	\item $\mathcal{PRK}_{6}43$: The 6-stage 4th-order splitting (\ref{prk1}), with the 3rd-order estimator (\ref{prk2}) and coefficients given by 
	(\ref{est.1}).
	\item $\mathcal{S}_{6}43$: The 6-stage 4th-order method-adjoint symmetric composition (\ref{eq:compint}) with coefficients (\ref{MetAdj1}) and 3rd-order estimator (\ref{MetAdj2}).
	\item $\mathcal{SS}_{5}43$: The 5-stage 4th-order symmetric composition (\ref{eq:compSS}) with coefficients (\ref{suz1}) and 3rd-order estimator (\ref{suz2}).
	\item $\mathcal{SS}_{11}65$: The 11-stage 6th-order symmetric composition (\ref{eq:compSS}) with coefficients (\ref{m65al}) and 5th-order 
	estimator (\ref{m65}).
	\item $\mathcal{SS}_{17}853$: The 17-stage 8th-order symmetric composition (\ref{eq:compSS}) with coefficients (\ref{krcl0})	with the 5th- and 3rd-order estimators (\ref{krcl1}) and (\ref{krcl2}).
\end{itemize}

These are compared with:
\begin{itemize}
	\item $\mathrm{eRKN}_{4}43$: the non-symmetric 4-stage 4th-order Runge--Kutta--Nystr{\"o}m (RKN) method with a 3rd-order estimator presented in \cite{calvo93tdo}. This method has an error estimator that is only valid for equations of the form (\ref{rkn1}), so that it cannot be used in particular for
	the Schr\"odinger equation. 
	\item $\mathrm{ePRK}_{5}43$: The 5-stage 4th-order splitting method given by the composition
\begin{equation}   \label{koch1}
  \psi_h = 
	\varphi_{b_5 h}^{[2]} \circ \varphi_{a_5 h}^{[1]} \circ   
	\varphi_{b_4 h}^{[2]} \circ \varphi_{a_4 h}^{[1]} \circ   
	\varphi_{b_3 h}^{[2]} \circ \varphi_{a_3 h}^{[1]} \circ   
	\varphi_{b_2 h}^{[2]} \circ \varphi_{a_2 h}^{[1]} \circ   
  \varphi_{b_1 h}^{[2]} \circ \varphi_{a_1 h}^{[1]}
\end{equation}
with the symmetry $b_{6-i}=a_i, \ i=1,2,\ldots,5$, and the 3rd-order estimator given by a similar composition sharing the first stages\footnote{The idea to consider estimators using a second composition sharing some of the stages was first proposed in \cite{koch13eeo}}
\begin{equation}   \label{koch2}
  \widetilde \psi_h = 
	\varphi_{\widetilde b_5 h}^{[2]} \circ \varphi_{\widetilde a_5 h}^{[1]} \circ   
	\varphi_{\widetilde b_4 h}^{[2]} \circ \varphi_{\widetilde a_4 h}^{[1]} \circ   
	\varphi_{\widetilde b_3 h}^{[2]} \circ \varphi_{\widetilde a_3 h}^{[1]} \circ   
	\varphi_{b_2 h}^{[2]} \circ \varphi_{a_2 h}^{[1]} \circ   
  \varphi_{b_1 h}^{[2]} \circ \varphi_{a_1 h}^{[1]}
\end{equation}
with $\widetilde a_i,\widetilde b_i, \ i=3,4,5$ chosen appropriately. The estimator requires three new evaluations. We take in particular the scheme\footnote{The corresponding coefficients are available at \url{http://www.asc.tuwien.ac.at/~winfried/splitting}}
\texttt{Emb 4/3 AK p}, in which case $\widetilde a_3=0$, so that only  two new evaluations are required and the overall cost is taken as
7 evaluations per step.
	\item RK6(5): the well known 8-stage Verner's method of order 6(5) (see Table 5.4 in \cite{hairer93sod}, page 181) that is implemented in the routine DVERK.
	\item DOP853: the 12-stage embedded Runge--Kutta method of order 8(5) by Dormand \& Prince \cite{dormand89prk} and improved  as the routine DOP853 in \cite{hairer93sod}.
\end{itemize}


%
Specifically, we consider as a test bench the two-dimensional
Kepler problem with Hamiltonian
\begin{equation}   \label{eq.HamKepler}
   H(q,p) = T(p) + V(q) = 
	\frac{1}{2} p^T p - \mu \frac{1}{r}.
\end{equation}
Here $q=(q_1,q_2), p=(p_1,p_2)$, $\mu=GM$, $G$ is the gravitational constant and $M$ is the sum of the masses of the two bodies. Taking $\mu=1$ and initial conditions
\begin{equation}\label{eq.1.12}
  q_1(0) = 1- e, \quad q_2(0) = 0, \quad p_1(0) = 0, \quad p_2(0) = \sqrt{\frac{1+e}{1-e}},
\end{equation}
if $0 \le e < 1$, then the total energy is $H=H_0=-1/2$, the solution is periodic with period $2 \pi$, and the trajectory is an ellipse of eccentricity $e$. 

The performance of an embedded Runge--Kuta method depends on the performance of the high order method used to propagate the solution, but also on the accuracy of the lower order one as well as how the error estimator approaches the true error of the high order method. Some times the error estimator is much larger than the true error and the algorithm uses smaller time steps than necessary to reach a given accuracy.
Some  other times, however, this error can be considerably smaller than the true error (usually due to cancellations because the methods share internal stages) and the algorithm takes longer time steps than required which lead to undesirable large errors.

In this example we integrate with a constant time step and compute the maximum true error
\[
  {\cal E}_1 =\max_n \| x(t_n)-x_n\|
\]
and the maximum error estimator
\[
  {\cal E}_2 =\max_n \| \widetilde{x}_n-x_n\|.
\]
An efficient method should give ${\cal E}_2 \sim {\cal E}_1$, while being both as small as possible at a given computational cost.

The integration is carried out in the time interval $t\in[0,20]$ with a constant time step, and this integration is repeated for different values of the time step and for several values of the eccentricity, in particular for  $e=\frac15, \ \frac25, \ \frac35, \ \frac45$. This is done first for RK6(5) (or
DVERK subroutine) and
the composition scheme $\mathcal{SS}_{11}65$. 

\begin{figure}[ht]
\begin{center}
\includegraphics[scale=0.7]{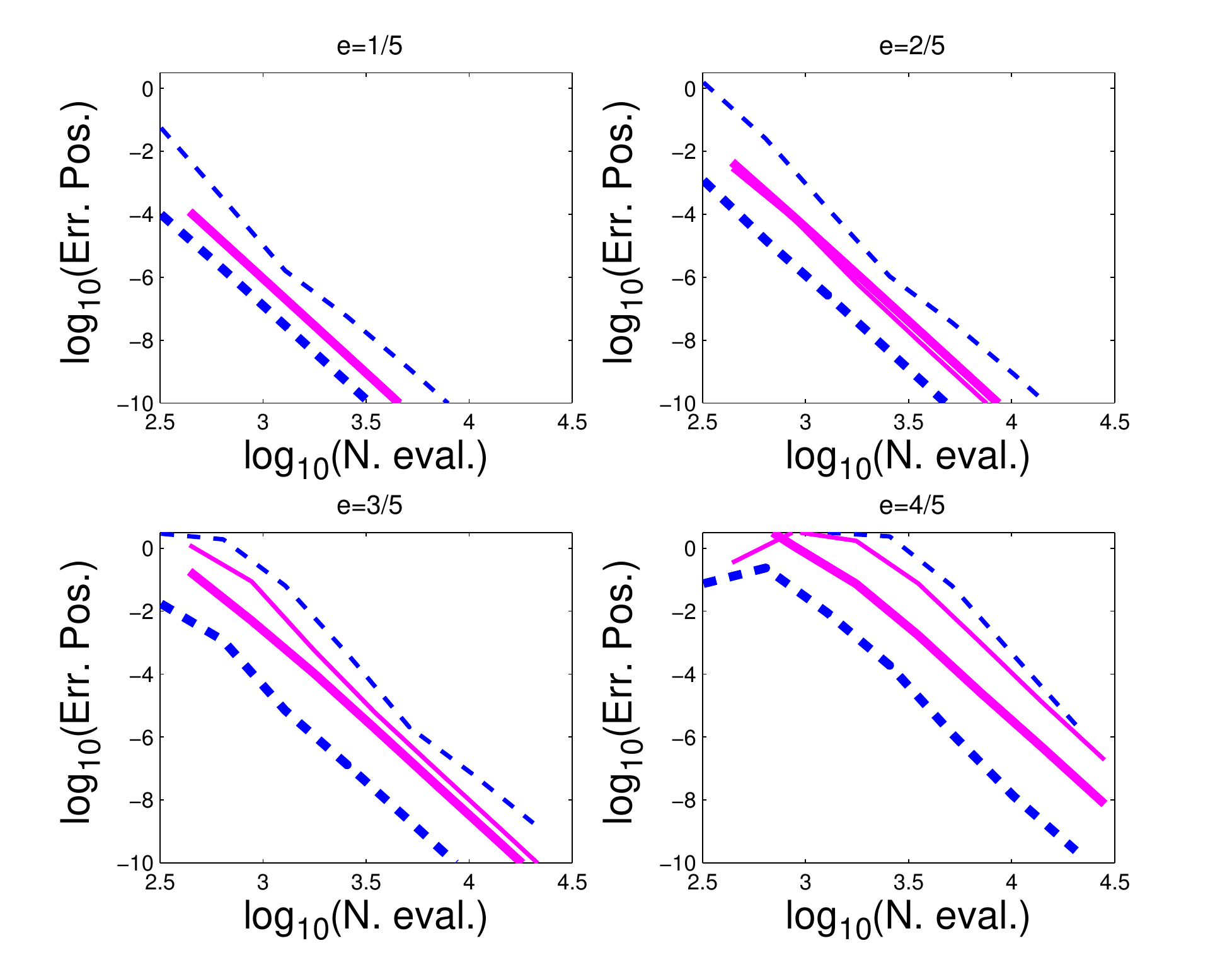}
\end{center}
\caption{{\bf Methods of order 6(5)}. Maximum error in positions, ${\cal E}_1$ (thin lines), and maximum error estimator, ${\cal E}_2$ (thick lines), versus the computational cost measured as the number of force evaluations in double logarithmic scale: (dashed lines) DVERK; and (solid lines) $\mathcal{SS}_{11}65$. 
}
 \label{fig1}
\end{figure}

Figure~\ref{fig1} shows in double logarithmic scale the error ${\cal E}_1$ (thin lines) and the estimate ${\cal E}_2$ (thick lines) versus
 the computational cost measured as the number of force evaluations. Dashed lines are obtained with RK6(5), whereas 
 solid lines correspond to $\mathcal{SS}_{11}65$.

We notice from the figure that the composition method is not only more accurate at the same cost (even for such a short time integration) but also the error estimator is much closer to the true error.
The error estimator of DVERK is very optimistic: ${\cal E}_2$ is much smaller that ${\cal E}_1$, especially when the eccentricity takes large values (and thus 
adjusting the step size is increasingly relevant). The reason lies in the fact that both $x_n$ and $\widetilde{x}_n$ are computed using very similar procedures, since they share the intermediate stages. 
This is not the case for the error estimators proposed here, and thus the error ${\cal E}_2$ is reasonably close to the true error of the method, even
when the coefficients for this specific method are not particularly small.

\begin{figure}[ht]
\begin{center}
\includegraphics[scale=0.75]{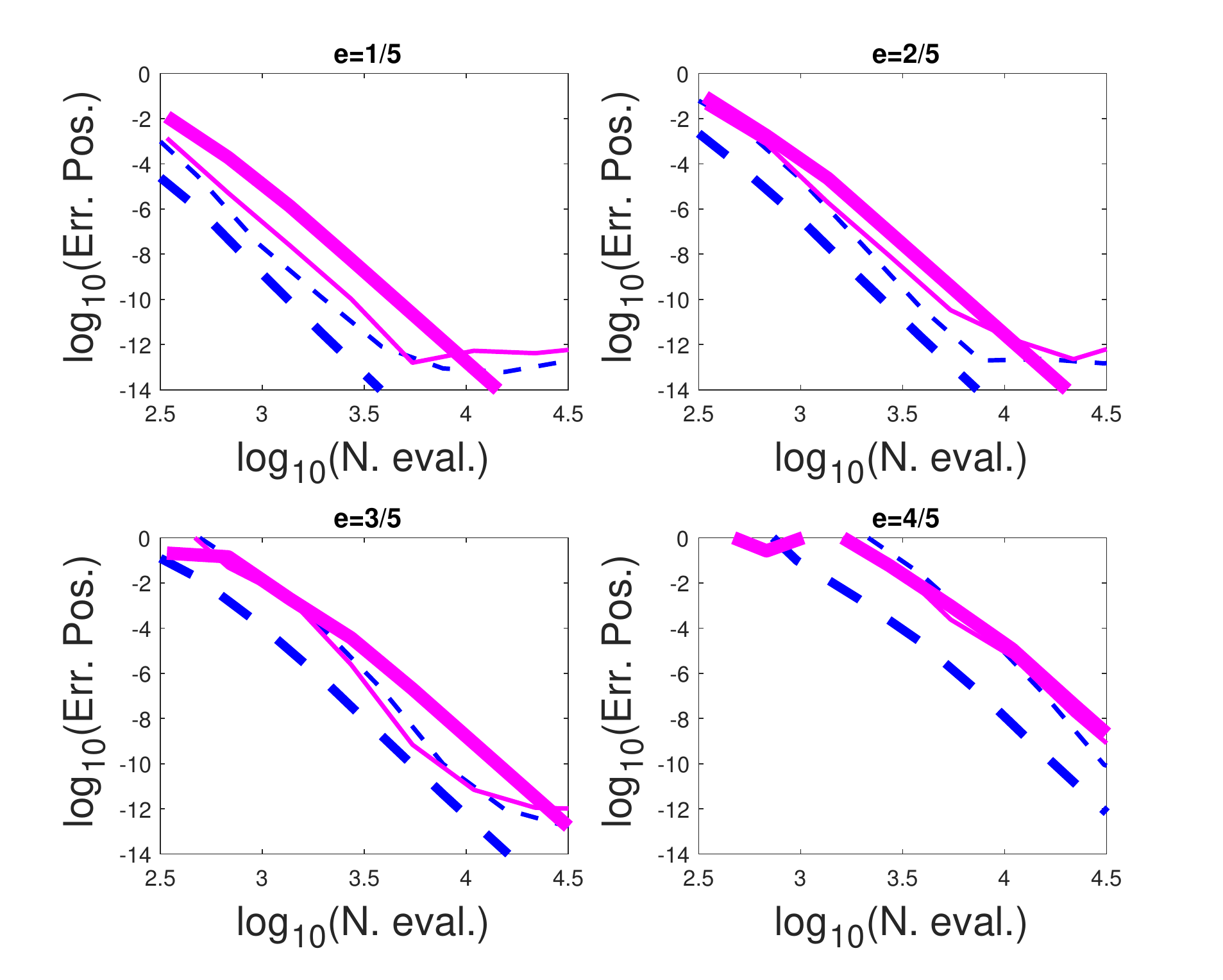}
\end{center}
\caption{{\bf Methods of order 8(5)(3)}. Maximum error in positions, ${\cal E}_1$ (thin lines), and maximum error estimator, ${\cal E}_2$ (thick lines), versus the computational cost measured as the number of force evaluations in double logarithmic scale: (dashed lines) DOP853; and (solid lines) $\mathcal{SS}_{17}853$. 
}
 \label{fig2}
\end{figure}

Next the same numerical experiment is carried out again, but this time with DOP853 and the composition scheme $\mathcal{SS}_{17}853$. Figure~\ref{fig2} shows the results obtained.

We observe that, for this example, the symplectic composition method is as efficient as the 8th-order RK method even for such a short time integration. In addition, our error estimator for the composition method is closer to the true error providing a better error estimator and as a result allowing to choose more appropriate time steps.

\begin{figure}[ht]
\begin{center}
\includegraphics[scale=0.75]{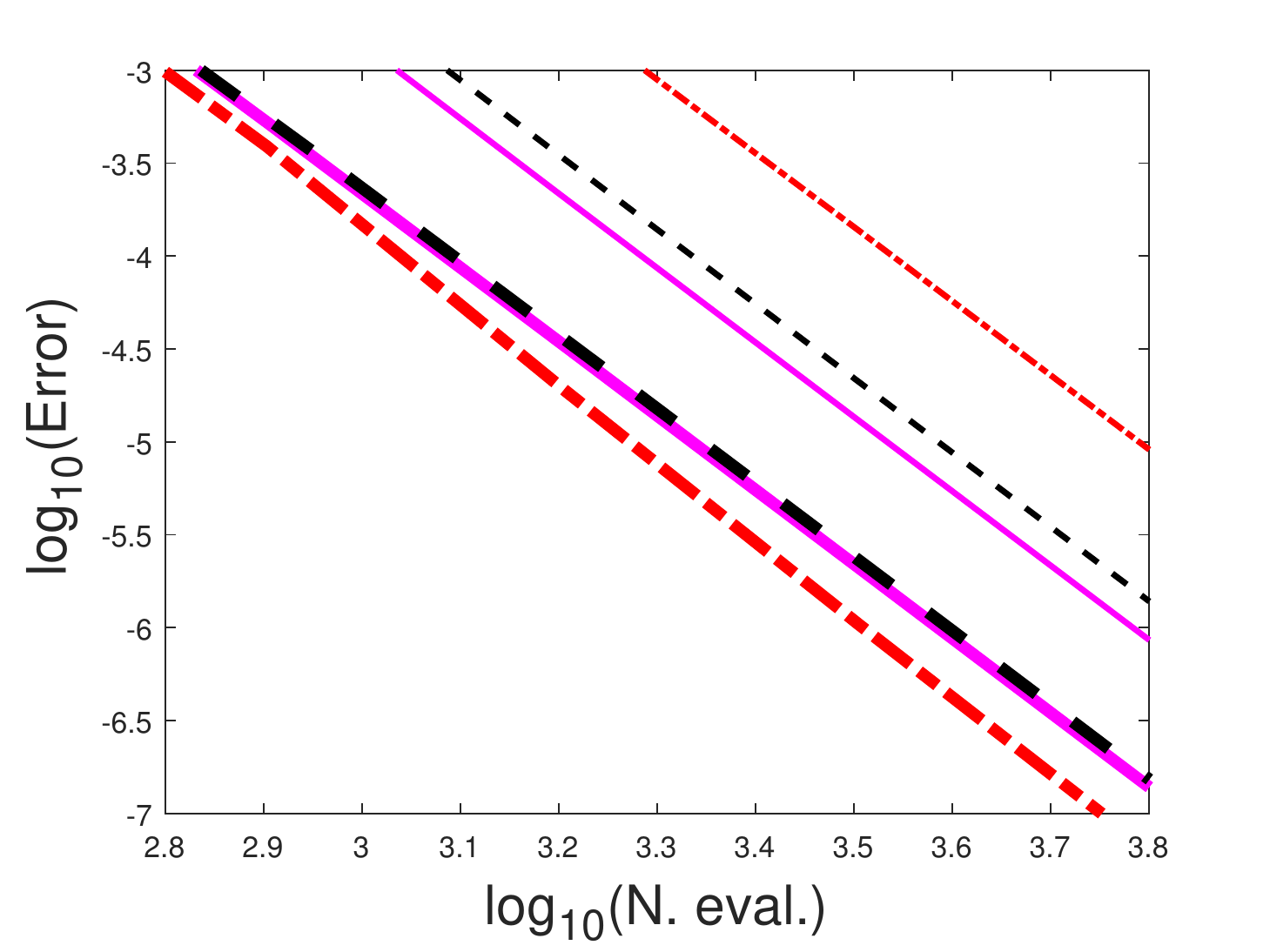}
\end{center}
\caption{{\bf General methods of order 4(3)}. Comparison of the true error (thin lines) and the estimator (thick lines)
for $e=1/2$ and the following schemes: 
$\mathcal{PRK}_{6}43$ (dashed lines); 
$\mathrm{ePRK}_{5}43$ (dot-dashed lines); and 
$\mathcal{SS}_{5}43$  (solid lines). 
}
 \label{fig3}
\end{figure}

Next we compare the results achieved by methods of order 4(3) that are valid for general splitting methods and symmetric-symmetric compositions.
This is shown in Figure~\ref{fig3} for eccentricity $e=1/2$ in eq. (\ref{eq.1.12}): 
 $\mathcal{PRK}_{6}43$ (dashed lines); 
$\mathrm{ePRK}_{5}43$ (dot-dashed lines); and 
$\mathcal{SS}_{5}43$  (solid lines). We observe that the embedded scheme $\mathrm{ePRK}_{5}43$ provides an exceedingly optimistic error estimator as well as a lower performance due to its higher cost per step.

\begin{figure}[ht]
\begin{center}
\includegraphics[scale=0.75]{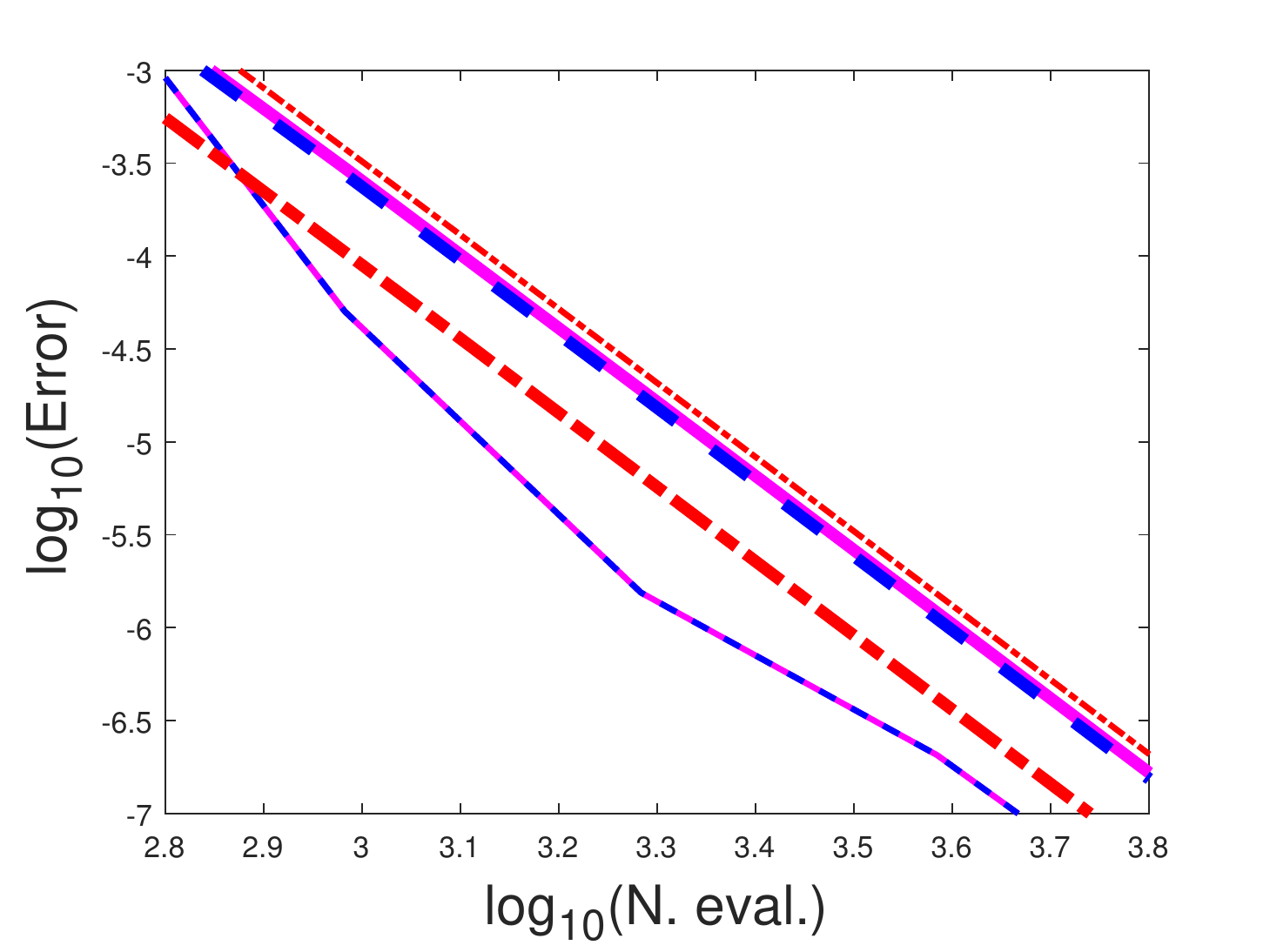}
\end{center}
\caption{{\bf RKN methods of order 4(3)}. Same as Figure~\ref{fig3} for the following methods: 
$\mathcal{S}_{6}43$ (dashed lines); 
$\mathrm{eRKN}_{4}43$ (dot-dashed lines);
 and 
$\mathcal{RKN}_{6}43$ (solid lines). 
}
 \label{fig4}
\end{figure}

Finally, Figure~\ref{fig4} shows the same results as Figure~\ref{fig3} for the RKN methods of order 4(3) and the composition method-adjoint obtained from the coefficients of the 6-stage RKN method and the relation (\ref{eq:alphaab}). It provides the same results for the 4th-order method, but different outputs for the estimator. Specifically, we collect the results obtained with 
$\mathcal{S}_{6}43$ (dashed lines), $\mathrm{eRKN}_{4}43$ (dot-dashed lines), and $\mathcal{RKN}_{6}43$ (solid lines). 
We observe that the scheme $\mathrm{eRKN}_{4}43$ provides an optimistic error estimator as well as a lower performance.

\section{Concluding remarks}

In this work we have proposed a procedure to estimate the local error of splitting and composition methods based on the construction of a 
second lower order integrator by linear combinations of the intermediate outputs of the original scheme. The difference can then be combined
with standard strategies of automatic step size control \cite{hairer93sod} to use the original splitting and composition methods with adaptive step size
along the integration. In contrast with other approaches, the proposed strategy does not increase the computational cost of the overall scheme and provides
a reliable estimate of the error, so that it can be safely used in problems where keeping the step size constant is not of paramount importance, 
such as it is the case in certain partial differential equations of evolution. In any event, in that case one should use a very precise discretization
in space to guarantee that the main source of error originates when integrating in time.

We should remark in particular the good properties exhibited by the estimator constructed for the 17-stage 8th-order composition scheme
(\ref{eq:compSS}) with coefficients (\ref{krcl0}) in comparison with the well known routine DOP853. Taking into account that even more
efficient composition methods involving 19 and 21 stages do exist within this class, we conclude that these can constitute a worthwhile
alternative for integrating problems when high accuracy is required.

The error estimator proposed here coupled with a variable step size strategy could be most useful for the application of splitting methods
for solving the Schr\"odinger eigenvalue problem with the imaginary time propagation technique, in order to reduce the overall computational
cost, as illustrated e.g. in \cite{bader13sts}.

Although only several representative schemes have been considered, it is clear that the same strategy can be applied to any other
splitting and composition method. In particular, we can also construct estimators for the high-order
methods with complex coefficients collected in \cite{blanes13oho} and schemes involving double commutators, such as those presented in
\cite{chin97sif,laskar01hos,omelyan03sai},
as long as they involve a sufficiently large number of intermediate stages to form
the required linear combinations.

\section*{Acknowledgements}
Part of this work was developed during a research stay at the \emph{Wolfgang Pauli Institute Vienna}; 
the authors are grateful to the director Norbert Mauser and the staff members for their support and hospitality.
The first two authors acknowledge funding by the Ministerio de Econom{\'i}a y Competitividad (Spain) through project MTM2016-77660-P (AEI/FEDER, UE).


\end{document}